\documentclass[12pt]{article}
\usepackage{latexsym,amssymb,amsmath,enumerate,geometry,float,cite,tikz,setspace}
\newtheorem{definition}{Definition}
\newtheorem{theorem}[definition]{Theorem}
\newtheorem{corollary}[definition]{Corollary}
\newtheorem{lemma}[definition]{Lemma}
\newtheorem{proposition}[definition]{Proposition}

\newtheorem{claim}{Claim}

\usepackage{lineno}

\begin{document}


\onehalfspace

\title{Averaging $2$-Rainbow Domination and Roman Domination}

\author{Jos\'{e} D. Alvarado$^1$, Simone Dantas$^1$, and Dieter Rautenbach$^2$}

\date{}

\maketitle

\begin{center}
{\small 
$^1$ Instituto de Matem\'{a}tica e Estat\'{i}stica, Universidade Federal Fluminense, Niter\'{o}i, Brazil,
\texttt{josealvarado.mat17@gmail.com, sdantas@im.uff.br}\\[3mm]
$^2$ Institute of Optimization and Operations Research, Ulm University, Ulm, Germany,
\texttt{dieter.rautenbach@uni-ulm.de}
}
\end{center}

\begin{abstract}
For a graph $G$, let $\gamma_{r2}(G)$ and $\gamma_R(G)$ denote the $2$-rainbow domination number and the Roman domination number, respectively. 
Fujita and Furuya (Difference between 2-rainbow domination and Roman domination in graphs, 
Discrete Applied Mathematics 161 (2013) 806-812) proved
$\gamma_{r2}(G)+\gamma_R(G)\leq \frac{6}{4}n(G)$
for a connected graph $G$ of order $n(G)$ at least $3$.
Furthermore, they conjectured
$\gamma_{r2}(G)+\gamma_R(G)\leq \frac{4}{3}n(G)$
for a connected graph $G$ of minimum degree at least $2$ that is distinct from $C_5$.
We characterize all extremal graphs for their inequality and prove their conjecture.
\end{abstract}

{\small 

\medskip


\medskip


}

\section{Introduction}\label{section1}

We consider finite, simple, and undirected graphs and use standard terminology and notation.

Rainbow domination of graphs was introduced in \cite{brhera}.
Here we consider the special case of $2$-rainbow domination. 
A {\it $2$-rainbow dominating function of a graph $G$} is a function $f:V(G)\to 2^{\{ 1,2\}}$ 
such that $\bigcup_{v\in N_G(u)}f(v)=\{ 1,2\}$ for every vertex $u$ of $G$ with $f(u)=\emptyset$.
The {\it weight of $f$} is $\sum_{u\in V(G)}|f(u)|$.
The {\it $2$-rainbow domination number $\gamma_{r2}(G)$ of $G$} 
is the minimum weight of a $2$-rainbow dominating function of $G$.
Roman domination was introduced in \cite{st}.
A {\it Roman dominating function of a graph $G$} is a function $g:V(G)\to \{ 0,1,2\}$ 
such that every vertex $u$ of $G$ with $g(u)=0$
has a neighbor $v$ with $g(v)=2$.
The {\it weight of $g$} is $\sum_{u\in V(G)}g(u)$.
The {\it Roman domination number $\gamma_R(G)$ of $G$} 
is the minimum weight of a Roman dominating function of $G$.

The definitions of the above two types of dominating functions have some obvious similarities;
vertices contribute either $0$ or $1$ or $2$ to the weight of these functions;
for vertices that contribute $0$, their neighbors contribute at least $2$ in total;
vertices that contribute $1$ do not impose any condition on their neighbors; and
vertices that contribute $2$ satisfy the requirements of all their neighbors that contribute $0$.
Nevertheless, while vertices that contribute $1$ are useless for their neighbors in Roman domination,
they can satisfy `half' the requirements of their neighbors in $2$-rainbow domination.

As observed in \cite{wuxi,chra} 
the two domination parameters are related by the following simple inequalities 
\begin{eqnarray}\label{e1}
\gamma_{r2}(G)\leq \gamma_R(G)\leq \frac{3}{2}\gamma_{r2}(G).
\end{eqnarray}
In fact, if $g$ is a Roman dominating function of a graph $G$ of weight $w$, then 
$$
f:V(G)\to 2^{\{ 1,2\}}:
u\mapsto
\left\{
\begin{array}{cl}
\emptyset & \mbox{, if $g(u)=0$,}\\
\{ 1\} & \mbox{, if $g(u)=1$, and}\\
\{ 1,2\} & \mbox{, if $g(u)=2$.}
\end{array}
\right.
$$
is a $2$-rainbow dominating function of $G$ of weight $w$, 
which implies $\gamma_{r2}(G)\leq \gamma_R(G)$.
Similarly,
if $f$ is a $2$-rainbow dominating function of $G$ of weight $w$
such that $|f^{-1}(\{ 1\})|\geq |f^{-1}(\{ 2\})|$, then 
$$
g:V(G)\to \{ 0,1,2\}:
u\mapsto
\left\{
\begin{array}{cl}
0 & \mbox{, if $f(u)=\emptyset$,}\\
1  & \mbox{, if $f(u)=\{ 1\}$, and}\\
2 & \mbox{, if $f(u)\in \{ \{ 2\},\{ 1,2\}\}$.}
\end{array}
\right.
$$
is a Roman dominating function of $G$ of weight at most $3w/2$, 
which implies $\gamma_R(G)\leq \frac{3}{2}\gamma_{r2}(G)$.

The following result summarizes known tight bounds for the two parameters \cite{wuxi,chkiprwe,fufu}. 

\begin{theorem}\label{theorembounds}
Let $G$ be a connected graph of order $n(G)$ at least $3$.
\begin{enumerate}[(i)]
\item $\gamma_{r2}(G)\leq \frac{3}{4}n(G)$ \cite{wuxi}.
\item $\gamma_R(G)\leq \frac{4}{5}n(G)$ \cite{chkiprwe}.
\item If $G$ has minimum degree at least $2$, then $\gamma_{r2}(G)\leq \frac{2}{3}n(G)$ \cite{fufu}.
\item If $G$ has order at least $9$ and minimum degree at least $2$, then $\gamma_R(G)\leq \frac{8}{11}n(G)$ \cite{chkiprwe}.
\end{enumerate} 
\end{theorem}
Also bounds on linear combinations of the parameters were considered.

\begin{theorem}[Fujita and Furuya \cite{fufu}]\label{theoremaverage1}
If $G$ is a connected graph of order $n(G)$ at least $3$, then $\gamma_{r2}(G)+\gamma_R(G)\leq \frac{6}{4}n(G)$.
\end{theorem}
In view of Theorem \ref{theorembounds}(i) and (ii),
one would expect an upper bound on $(\gamma_{r2}(G)+\gamma_R(G))/(2n(G))$
that is somewhere between $3/4$ and $4/5$.
Theorem \ref{theoremaverage1} is slightly surprising 
as it shows that this upper bound has the smallest possible value, namely $3/4$.
In fact, (\ref{e1}) and Theorem \ref{theoremaverage1} imply Theorem \ref{theorembounds}(i).

Our first result is the following.

\begin{theorem}\label{theorem1}
If $G$ is a connected graph of minimum degree at least $2$ that is distinct from $C_5$, 
then $\gamma_{r2}(G)+\gamma_R(G)\leq \frac{4}{3}n(G)$.
\end{theorem}
Theorem \ref{theorem1} confirms a conjecture of Fujita and Furuya (Conjecture 2.11 in \cite{fufu}).
Similarly as for Theorem \ref{theoremaverage1},
it is again slightly surprising that the upper bound on $(\gamma_{r2}(G)+\gamma_R(G))/(2n(G))$
in Theorem \ref{theorem1} has the smallest of the possible values suggested by Theorem \ref{theorembounds}(iii) and (iv),
namely $2/3$.
Note that (\ref{e1}) and Theorem \ref{theorem1} imply Theorem \ref{theorembounds}(iii).

Another result concerning linear combinations of the parameters is the following.
\begin{proposition}[Chellali and Rad \cite{chra}]\label{theoremaverage2}
There is no constant $c$ such that $2\gamma_{r2}(G)+\gamma_R(G)\leq 2n(G)+c$ for every connected graph $G$.
\end{proposition}
Chellali and Rad posed the problem (Problem 13 in \cite{chra}) to find a sharp upper bound on
$2\gamma_{r2}(G)+\gamma_R(G)$ for connected graphs $G$ of order at least $3$.
In fact, (\ref{e1}) and Theorem \ref{theoremaverage1} immediately imply the following.

\begin{corollary}\label{corollary1}
If $G$ is a connected graph of order $n(G)$ at least $3$, 
then $2\gamma_{r2}(G)+\gamma_R(G)\leq \frac{9}{4}n(G)$.
\end{corollary}
Corollary \ref{corollary1} is sharp and hence solves the problem posed by Chellali and Rad.

As our second result we characterize all extremal graphs for Theorem \ref{theoremaverage1},
all of which are also extremal for Corollary \ref{corollary1}.

\section{Results and Proofs}

Our proof of Theorem \ref{theorem1} relies on an elegant approach from \cite{chkiprwe}.
We also use the reductions described in Lemma 4.1 in \cite{chkiprwe}.
Unfortunately, the proofs of (b) and (c) of Lemma 4.1 in \cite{chkiprwe} are not completely correct;
the graphs $G'$ considered in these proofs may have vertices of degree less than $2$. 
We incorporate corrected proofs for these reductions as claims within the proof of Theorem \ref{theorem1}.

\begin{lemma}\label{lemma1}
Let $G$ be a graph that contains an induced path $P$ of order $5$ whose internal vertices have degree $2$.
If $G'$ arises from $G$ by contracting three edges of $P$, then 
$n(G')=n(G)-3$, 
$\gamma_{r2}(G)\leq \gamma_{r2}(G')+2$, and
$\gamma_R(G)\leq \gamma_R(G')+2$.
\end{lemma}
{\it Proof:} Let $P:xuvwy$,
that is, $G'$ arises from $G$ by deleting $u$, $v$, and $w$, and adding the edge $xy$.
Clearly, $n(G')=n(G)-3$.

Let $f$ be a $2$-rainbow dominating function of $G'$.
If $f(x),f(y)\not=\emptyset$ or $f(x)=f(y)=\emptyset$, 
then setting $f(u)=f(w)=\emptyset$ and $f(v)=\{ 1,2\}$ 
extends $f$ to a $2$-rainbow dominating function of $G$.
Now we assume that $1\in f(x)$ and $f(y)=\emptyset$.
If $f(x)=\{ 1\}$, 
then setting $f(u)=\emptyset$, $f(v)=\{ 2\}$, and $f(w)=\{ 1\}$
extends $f$ to a $2$-rainbow dominating function of $G$.
Finally, if $f(x)=\{ 1,2\}$, 
then setting $f(u)=f(v)=\emptyset$ and $f(w)=\{ 1,2\}$
extends $f$ to a $2$-rainbow dominating function of $G$.
By symmetry, this implies $\gamma_{r2}(G)\leq \gamma_{r2}(G')+2$.

Let $g$ be a Roman dominating function of $G'$.
If $g(x)=2$ and $g(y)=0$, 
then setting $g(u)=g(v)=0$ and $g(w)=2$ 
extends $g$ to a Roman dominating function of $G$.
Now we assume that $\{ g(x),g(y)\}\not=\{ 0,2\}$.
Setting $g(u)=g(w)=0$ and $g(v)=2$ 
extends $g$ to a Roman dominating function of $G$.
By symmetry, this implies $\gamma_R(G)\leq \gamma_R(G')+2$.
$\Box$

\medskip

\noindent A {\it spider} is a graph that arises 
by iteratively subdividing the edges of a star of order at least $4$ arbitrarily often.
The {\it center} of a spider is its unique vertex of degree at least $3$.
The {\it legs} of a spider are the maximal paths starting at its center.
A leg is {\it good} if its length is not a multiple of $3$.

\begin{lemma}\label{lemma3}
If $G$ is a spider with at least three good legs,
then $\gamma_{r2}(G)+\gamma_R(G)\leq \frac{4}{3}n(G)$.
\end{lemma}
{\it Proof:}
If some leg of $G$ has length at least $3$, 
then contracting three edges of this leg yields a spider $G'$ 
with at least three good legs such that
$n(G')=n(G)-3$, and $\gamma_{r2}(G)+\gamma_R(G)\leq \gamma_{r2}(G')+\gamma_R(G')+4$.
By an inductive argument, we may therefore assume that all legs of $G$ have length $1$ or $2$.
Let $G$ have exactly $\ell_1$ legs of length $1$ and exactly $\ell_2$ legs of length $2$.
By the hypothesis, $\ell_1+\ell_2\geq 3$. 

If $\ell_1=0$, then
$\gamma_{r2}(G)+\gamma_R(G)\leq (1+\ell_2)+(2+\ell_2)=3+2\ell_2$
and $n(G)=1+2\ell_2$.
Since $\ell_2\geq 3$, we obtain $\frac{3+2\ell_2}{1+2\ell_2}\leq \frac{9}{7}<\frac{4}{3}$,
and hence $\gamma_{r2}(G)+\gamma_R(G)\leq \frac{4}{3}n(G)$.
Now let $\ell_1\geq 1$.
We have 
$\gamma_{r2}(G)+\gamma_R(G)\leq (2+\ell_2)+(2+\ell_2)=4+2\ell_2$
and $n(G)=1+\ell_1+2\ell_2$.
Since $\frac{4+2\ell_2}{1+\ell_1+2\ell_2}<\frac{4+2\ell_2}{1+(\ell_1-1)+2\ell_2}$,
we may assume that either $\ell_1=1$ 
or $\ell_1\geq 2$ and $\ell_1+\ell_2=3$.

If $\ell_1=1$, then $\ell_2\geq 2$, and we obtain 
$\frac{4+2\ell_2}{1+\ell_1+2\ell_2}=\frac{4+2\ell_2}{2+2\ell_2}\leq \frac{4}{3}$.
If $\ell_1=2$, then $\ell_2=1$, and we obtain 
$\frac{4+2\ell_2}{1+\ell_1+2\ell_2}=\frac{4}{5}$.
If $\ell_1=3$, then $\ell_2=0$, and we obtain 
$\frac{4+2\ell_2}{1+\ell_1+2\ell_2}=\frac{4}{4}$.
Therefore, in all these cases, $\gamma_{r2}(G)+\gamma_R(G)\leq \frac{4}{3}n(G)$.
$\Box$

\medskip

\noindent We proceed to the proof of our first result.

\medskip

\noindent {\it Proof of Theorem \ref{theorem1}:} Let $G$ be a counterexample of minimum order.
A {\it branch vertex} is a vertex of degree at least $3$.
A {\it thread} is 
either a path between two branch vertices whose internal vertices have degree $2$
or a cycle with exactly one branch vertex.
By Lemma \ref{lemma1}, 
every thread that is a path has length at most $4$, and  
every thread that is a cycle has length at most $5$.

Since 
$\gamma_{r2}(C_3)+\gamma_R(C_3)=2+2\leq \frac{4}{3}\cdot 3$,
$\gamma_{r2}(C_4)+\gamma_R(C_4)=2+3\leq \frac{4}{3}\cdot 4$,
$\gamma_{r2}(C_8)+\gamma_R(C_8)=4+6\leq \frac{4}{3}\cdot 8$,
and, by Lemma \ref{lemma1},
$\gamma_{r2}(C_n)+\gamma_R(C_n)\leq \gamma_{r2}(C_{n-3})+\gamma_R(C_{n-3})+4$
for $n\geq 6$,
a simple inductive argument implies $\gamma_{r2}(C_n)+\gamma_R(C_n)\leq \frac{4}{3}n$
for every $n\geq 3$ that is distinct from $5$.
Hence, the maximum degree of $G$ is at least $3$.

If $G$ arises from $C_5$ by adding at least one edge, 
then $\gamma_{r2}(G)+\gamma_R(G)\leq 3+3<\frac{4}{3}n(G)$.
Hence, we may assume that $G$ does not have this structure.
Since removing edges can not decrease any of the two domination parameters,
we may assume that every edge between two branch vertices is a bridge,
that is, an edge whose removal increases the number of components.
Therefore, again by Lemma \ref{lemma1}, 
no thread that is a path has length $4$,
that is, every thread that is a path has length $1$, $2$, or $3$.

If $u$ and $v$ are adjacent branch vertices and none of the two components $G_1$ and $G_2$ of $G-uv$ is $C_5$,
then the choice of $G$ implies 
$\gamma_{r2}(G)+\gamma_R(G)\leq
(\gamma_{r2}(G_1)+\gamma_R(G_1))+(\gamma_{r2}(G_2)+\gamma_R(G_2))
\leq \frac{4}{3}n(G_1)+\frac{4}{3}n(G_2)=\frac{4}{3}n(G)$.
Hence, for every edge $uv$ between branch vertices $u$ and $v$, $G-uv$ contains a component that is $C_5$.

\begin{claim}\label{claim1}
No two branch vertices are joined by two threads of length $2$.
\end{claim}
{\it Proof of Claim \ref{claim1}:}
For a contradiction, we assume that the two branch vertices $x$ and $y$
are joined by two threads of length $2$,
that is, $x$ and $y$ have at least two common neighbors of degree $2$.
Note that $x$ and $y$ are not adjacent and do not have a common neighbor of degree at least $3$.
Let $G'$ arise from $G$ by contracting all edges incident with the common neighbors of $x$ and $y$
to form a new vertex $z$. Since no thread of $G$ is a path of length $5$, $G'$ is not $C_5$.

First we assume that $N_G(x)\setminus N_G(y),N_G(y)\setminus N_G(x)\not=\emptyset$.
In this case, $\delta(G')\geq 2$ and $n(G')\leq n(G)-3$.
Let $f$ be a $2$-rainbow dominating function of $G'$.
If $f(z)=\emptyset$, then setting $f(x)=\{ 1\}$, $f(y)=\{ 2\}$, and $f(u)=\emptyset$ for $u\in N_G(x)\cap N_G(y)$
extends $f$ to a $2$-rainbow dominating function of $G$.
If $f(z)\not=\emptyset$, then 
setting $f(x)=f(z)$, $f(y)=\{ 1,2\}$, and $f(u)=\emptyset$ for $u\in N_G(x)\cap N_G(y)$
extends $f$ to a $2$-rainbow dominating function of $G$.
In both cases, we obtain $\gamma_{r2}(G)\leq \gamma_{r2}(G')+2$.
Similarly, we obtain $\gamma_R(G)\leq \gamma_R(G')+2$.
Therefore, by the choice of $G$,
$\gamma_{r2}(G)+\gamma_R(G)\leq 
\gamma_{r2}(G')+\gamma_R(G')+4
\leq \frac{4}{3}n(G')+4\leq \frac{4}{3}n(G)$.

Next we assume that $N_G(x)\setminus N_G(y)=\emptyset$.
This implies that $x$ and $y$ have at least three common neighbors.
Let $G''$ arise from $G'$ by adding two new vertices $z'$ and $z''$,
and adding the three new edges $zz'$, $z'z''$, and $z''z$.
Clearly, $\delta(G'')\geq 2$, $n(G'')\leq n(G)-2$, and $G''$ is not $C_5$.
Let $f$ be a $2$-rainbow dominating function of $G''$.
Setting 
$f(x)=\{ 1\}$,
$f(y)=\{ 1,2\}$,
and
$f(u)=\emptyset$ for $u\in N_G(x)\cap N_G(y)$
extends $f$ to a $2$-rainbow dominating function of $G$.
Since $|f(z)|+|f(z')|+|f(z'')|\geq 2$, this implies $\gamma_{r2}(G)\leq \gamma_{r2}(G'')+1$.
Similarly, we obtain $\gamma_R(G)\leq \gamma_R(G'')+1$.
Therefore, by the choice of $G$,
$\gamma_{r2}(G)+\gamma_R(G)\leq 
\gamma_{r2}(G'')+\gamma_R(G'')+2
\leq \frac{4}{3}n(G'')+2<\frac{4}{3}n(G)$.
$\Box$

\begin{claim}\label{claim2}
No two branch vertices are joined by two threads of length $3$.
\end{claim}
{\it Proof of Claim \ref{claim2}:}
For a contradiction, we assume that the two branch vertices $x$ and $y$
are joined by two threads of length $3$.
Note that $x$ and $y$ are not adjacent and do not have a common neighbor of degree at least $3$,
that is, every common neighbor of $x$ and $y$ is the internal vertex of a thread of length $2$ between $x$ and $y$.
Let $G'$ arise from $G$ by 
contracting the edges of all threads between $x$ and $y$
to form a new vertex $z$,
adding two new vertices $z'$ and $z''$,
and 
adding the three new edges $zz'$, $z'z''$, and $z''z$.
Clearly, $\delta(G')\geq 2$, $n(G')\leq n(G)-3$, and $G'$ is not $C_5$.

Let $f$ be a $2$-rainbow dominating function of $G'$.
Setting 
$f(x)=f(y)=\{ 1,2\}$
and
$f(u)=\emptyset$ for $u\in V(C)\setminus \{ x,y\}$
extends $f$ to a $2$-rainbow dominating function of $G$.
Since $|f(z)|+|f(z')|+|f(z'')|\geq 2$,
this implies $\gamma_{r2}(G)\leq \gamma_{r2}(G')+2$.
Similarly, we obtain $\gamma_R(G)\leq \gamma_R(G')+2$.
Therefore, by the choice of $G$,
$\gamma_{r2}(G)+\gamma_R(G)\leq 
\gamma_{r2}(G')+\gamma_R(G')+4
\leq \frac{4}{3}n(G')+4\leq \frac{4}{3}n(G)$.
$\Box$

\begin{claim}\label{claim3}
$G$ does not have a thread $C$ that is a cycle of length $5$
whose branch vertex $u$ has degree exactly $3$
and is joined by a thread $P$ of length $2$ or $3$ to another branch vertex $v$.
\end{claim}
{\it Proof of Claim \ref{claim3}:}
Let $G'=G[(V(C)\cup V(P))\setminus \{ v\}]$ and $G''=G-V(G')$.
Since $G'$ contains a spanning subgraph that is a spider with three good legs, 
Lemma \ref{lemma3} implies $\gamma_{r2}(G')+\gamma_R(G')\leq \frac{4}{3}|V(G')|$.
If $G''$ is not $C_5$, then, by the choice of $G$,
$\gamma_{r2}(G)+\gamma_R(G)\leq 
(\gamma_{r2}(G')+\gamma_R(G'))
+
(\gamma_{r2}(G'')+\gamma_R(G''))
\leq \frac{4}{3}|V(G')|+\frac{4}{3}(n(G)-|V(G')|)=\frac{4}{3}n(G)$.
If $G''$ is $C_5$ and $P$ has length $2$, then
$\gamma_{r2}(G)+\gamma_R(G)=6+8<\frac{4}{3}\cdot 11=\frac{4}{3}n(G)$.
Finally,
if $G''$ is $C_5$ and $P$ has length $3$, then
$\gamma_{r2}(G)+\gamma_R(G)=8+8\leq \frac{4}{3}\cdot 12=\frac{4}{3}n(G)$.
$\Box$

\begin{claim}\label{claim4}
There is a set $E$ of edges such that $G-E$ has one component $S$ that is a spider with only good legs
and all remaining components of $G-E$ have minimum degree at least $2$ and are distinct from $C_5$.
\end{claim}
{\it Proof of Claim \ref{claim4}:}
If $G$ consists of two copies of $C_5$ together with one bridge,
then $\gamma_{r2}(G)+\gamma_R(G)=6+7<\frac{4}{3}\cdot 10=\frac{4}{3}n(G)$.
Hence, $G$ does not have this structure.
Together with Claim \ref{claim3}, this implies that $G$ has some branch vertex $c$
such that either $c$ does not lie on a thread that is a cycle of length $5$ 
or $c$ has degree at least $4$.
We will describe the construction of $E$ starting with the empty set 
such that $c$ is the center of the spider $S$
and $S$ has $d_G(c)$ good legs,
that is, every edge of $G$ incident with $c$ will be the initial edge of a leg of $S$.
\begin{itemize}
\item For every thread $C$ of $G$ that is a cycle and contains $c$, 
add to $E$ exactly one edge of $C$ at maximum distance from $c$.
Since $C$ has length $3$, $4$, or $5$, this leads to two good legs for $S$.
\item For every thread $P$ of $G$ that is a path of length $1$ and contains $c$,
the only edge $e$ of $P$ is a bridge and the component $K$ of $G-e$ that contains the neighbor of $c$ is $C_5$.
Add to $E$ exactly one edge of $K$ that is incident with a neighbor of $c$.
This leads to one good leg of length $5$ for $S$.
\end{itemize}
Let $U$ be the set of branch vertices of $G$ 
that are joined to $c$ by at least one thread that is a path of length at least $2$.
Note that $U$ contains no neighbor of $c$.
By Claims \ref{claim1} and \ref{claim2}, 
no vertex in $U$ is joined to $c$ by more than two threads,
and
if some vertex $U$ is joined to $c$ by two threads, 
then these two threads have length $2$ and $3$, respectively.
For $u$ in $U$, let $E(u)$ denote the edges incident with $u$ that lie on threads between $c$ and $u$.
For a branch vertex $u$ distinct from $c$, let $p(u)$ be the number of threads between $c$ and $u$.
Note that for $u\in U$, we have $p(u)=|E(u)|$.

\begin{itemize}
\item Let $u_1$ be a vertex of degree exactly $3$ with $p(u_1)=2$.

Let $Q_1$ denote the unique third thread starting at $u_1$ and 
joining $u_1$ with some vertex $u_2$ distinct from $c$.
Let $u_1,u_2,\ldots,u_k$ be a maximal sequence of distinct branch vertices such that
for $2\leq i\leq k-1$,
the vertex $u_i$ has degree exactly $p(u_i)+2$ and
$u_i$ is joined to $u_{i+1}$ by exactly one thread $Q_i$.
Since the sequence $u_1,u_2$ with only two elements satisfies all requirements trivially, 
such a maximal sequence is well defined.
We call such a maximal sequence {\it special}.

Note that for $2\leq i\leq k-1$, the vertex $u_i$ belongs to $U$, 
and is joined to $u_{i-1}$ by the thread $Q_{i-1}$
and to $u_{i+1}$ by the thread $Q_i$.
Furthermore, note that $u_k$ not necessarily belongs to $U$.

By the maximality of the sequence, 
the vertex $u_k$ has either degree $p(u_k)+1$
or degree at least $p(u_k)+3$. 
\begin{itemize}
\item First we assume that $d_G(u_k)=p(u_k)+1$.

Add to $E$ the set $E(u_i)$ for every $2\leq i\leq k$.
Since the two threads between $c$ and $u_1$ have length $2$ and $3$,
adding to $E$ one of the two edges in $E(u_1)$ leads to $p(u_1)+\cdots +p(u_k)$ good legs for $S$.
\item Next we assume that $d_G(u_k)\geq p(u_k)+3$.

If $u_k$ has degree exactly $p(u_k)+3$ and belongs to a thread $C$ that is a cycle of length $5$,
then add to $E$ 
\begin{itemize}
\item $E(u_i)$ for every $2\leq i\leq k$, 
\item one edge of $C$ incident with $u_k$, and 
\item one of the two edges in $E(u_1)$.
\end{itemize}
Again the edge from $E(u_1)$ can be chosen such that we obtain
$p(u_1)+\cdots +p(u_k)$ good legs for $S$.

If $u_k$ does not have degree exactly $p(u_k)+3$ or does not belong to a thread $C$ that is a cycle of length $5$,
then add to $E$ 
\begin{itemize}
\item $E(u_i)$ for every $2\leq i\leq k$, 
\item the edge of $Q_{k-1}$ incident with $u_k$, and 
\item one of the two edges in $E(u_1)$.
\end{itemize}
Again the edge from $E(u_1)$ can be chosen such that we obtain
$p(u_1)+\cdots +p(u_k)$ good legs for $S$.
\end{itemize}
\end{itemize}
It remains to consider the vertices in $U$ that do not lie in some special sequence.
Note that each such vertex $u$ satisfies 
either $d_G(u)\geq 4$ 
or $d_G(u)=3$ and $p(u)=1$.
\begin{itemize}
\item Let $u$ in $U$ be such that $d_G(u)\geq 4$, $p(u)=1$, 
and $u$ does not belong to some special sequence.

Add to $E$ the unique edge in $E(u)$.
This leads to one good leg for $S$.
\item Let $u$ in $U$ be such that $d_G(u)\geq 4$, $p(u)=2$, 
and $u$ does not belong to some special sequence.

If $d_G(u)=4$ and $u$ belongs to a thread $C$ that is a cycle of length $5$,
then add to $E$ one edge of $C$ incident with $u$
and one of the two edges in $E(u)$.
Again, this last edge can be chosen such that we obtain two good legs for $S$.

If $d_G(u)>4$ or $u$ does not belong to a thread $C$ that is a cycle of length $5$,
then add to $E$ the set $E(u)$.
This leads to two good legs for $S$.
\item Finally, let $u$ in $U$ be such that $d_G(u)=3$, $p(u)=1$, 
and $u$ does not belong to some special sequence.

By Claim \ref{claim3}, $u$ does not belong to a thread $C$ that is a cycle of length $5$.
Add to $E$ the unique edge in $E(u)$.
This leads to one good leg for $S$.
\end{itemize}
This completes the construction of $E$. 
As argued above, the choice of $c$ implies that we obtain one good leg for $S$ for every edge of $G$ incident with $c$.
Furthermore, the construction of $E$ easily implies that $G-V(S)$ has the desired properties.
$\Box$

\medskip

\noindent Let $S$ be in Claim \ref{claim4} and let $R=G-V(S)$.
By the choice of $G$ and Lemma \ref{lemma3},
$\gamma_{r2}(G)+\gamma_R(G)\leq 
(\gamma_{r2}(S)+\gamma_R(S))+(\gamma_{r2}(R)+\gamma_R(R))
\leq \frac{4}{3}n(S)+\frac{4}{3}n(R)
=\frac{4}{3}n(G)$,
which completes the proof. $\Box$

\medskip

\noindent We proceed to our second result, the characterization of all extremal graphs for Theorem \ref{theoremaverage1}.
First we characterize the extremal trees.

For $k\in \mathbb{N}$, let ${\cal T}_k$ be the set of all trees $T$ that arise 
from $k$ disjoint copies $a_1b_1c_1d_1,\ldots,a_kb_kc_kd_k$ of the path of order $4$
by adding some edges between vertices in $\{ b_1,\ldots,b_k\}$,
that is, $T[\{ b_1,\ldots,b_k\}]$ is a tree.
Let ${\cal T}=\bigcup_{k\in \mathbb{N}}{\cal T}_k$.

\begin{lemma}\label{lemma4}
For $k\in \mathbb{N}$, $\gamma_{r2}(T)=\gamma_R(T)=3k$ for every $T\in {\cal T}_k$.
\end{lemma}
{\it Proof:} Let $T\in {\cal T}_k$.
We denote the vertices of $T$ as in the definition of ${\cal T}_k$.

Since 
$$
f_k:V(T)\to 2^{\{ 1,2\}}:
u\mapsto
\left\{
\begin{array}{cl}
\emptyset & \mbox{, if $u\in \{ a_1,\ldots,a_k\}\cup \{ c_1,\ldots,c_k\}$,}\\
\{ 1,2\} & \mbox{, if $u\in \{ b_1,\ldots,b_k\}$, and}\\
\{ 1\} & \mbox{, if $u\in \{ d_1,\ldots,d_k\}$}
\end{array}
\right.
$$
and
$$
g_k:V(T)\to \{ 0,1,2\}:
u\mapsto
\left\{
\begin{array}{cl}
0 & \mbox{, if $u\in \{ a_1,\ldots,a_k\}\cup \{ c_1,\ldots,c_k\}$,}\\
2 & \mbox{, if $u\in \{ b_1,\ldots,b_k\}$, and}\\
1 & \mbox{, if $u\in \{ d_1,\ldots,d_k\}$}
\end{array}
\right.
$$
are a $2$-rainbow dominating function and a Roman dominating function of $T$, respectively,
we obtain $\gamma_{r2}(T),\gamma_R(T)\leq 3k$.

For every $2$-rainbow dominating function $f$ and every Roman dominating function $g$ of $T$,
and for every $1\leq i\leq k$, it is easy to see that 
$|f(a_i)|+|f(b_i)|+|f(c_i)|+|f(d_i)|\geq 3$
and 
$g(a_i)+g(b_i)+g(c_i)+g(d_i)\geq 3$.
Hence, we obtain $\gamma_{r2}(T),\gamma_R(T)\geq 3k$.
$\Box$

\medskip

\noindent The following result is a strengthened version of Theorem 2.8 in \cite{fufu}.

\begin{lemma}\label{lemma5}
If $T$ is a tree of order $n(T)$ at least $3$, 
then $\gamma_{r2}(T)+\gamma_R(T)\leq \frac{6}{4}n(T)$
with equality if and only if $T\in {\cal T}$.
\end{lemma}
{\it Proof:} We prove the statement by induction on the order of $T$.
If $T$ is a path of order at most $5$, then the statement is easily verified.
If $T$ is a spider with at least three good legs, 
then Lemma \ref{lemma3} implies $\gamma_{r2}(T)+\gamma_R(T)\leq \frac{4}{3}n(T)< \frac{6}{4}n(T)$.
Hence, we may assume that $T$ is not any of these trees.

Let $u_1u_2\ldots u_\ell$ be a longest path in $T$ such that $d_G(u_2)$ is maximum possible.
Clearly, $\ell\geq 4$.

If $d_T(u_2)\geq 3$, then let $T'=T-(N_G[u_2]\setminus \{ u_3\})$.
Clearly, $\gamma_{r2}(T)\leq \gamma_{r2}(T')+2$ and $\gamma_R(T)\leq \gamma_R(T')+2$.
Since $T$ is not a spider with at least three good legs, $T'$ has order at least $3$.
Hence, by induction,  
$\gamma_{r2}(T)+\gamma_R(T)
\leq 4+\gamma_{r2}(T')+\gamma_R(T')
\leq 4+\frac{6}{4}n(T')
\leq 4+\frac{6}{4}(n(T)-3)
<\frac{6}{4}n(T)$.
Hence, we may assume that $d_T(u_2)=2$.

Let $T_1$ and $T_2$ be the components of $T-u_3u_4$ such that $T_1$ contains $u_3$.
Clearly, $T_1$ is either a path of order between $3$ and $5$,
or a spider with at least three good legs.
Since $T$ is not a path of order at most $5$ or a spider with at least three good legs,
$T_2$ has order at least $3$. 
If $T_1$ is not a path of order $4$, then $T\not\in {\cal T}$ and, by induction,
$\gamma_{r2}(T)+\gamma_R(T)
\leq (\gamma_{r2}(T_1)+\gamma_R(T_1))+(\gamma_{r2}(T_2)+\gamma_R(T_2))
<\frac{6}{4}n(T_1)+\frac{6}{4}n(T_2)
=\frac{6}{4}n(T)$.
Hence, we may assume that $T_1$ is a path of order $4$.
By a similar argument we obtain $T_2\in {\cal T}$.
Let $T_2\in {\cal T}_{k-1}$ and denote the vertices of $T_2$ as in the definition of ${\cal T}_{k-1}$.
Let $T_1$ be $a_kb_kc_kd_k$ such that $u_3=b_k$.
Let $a_ib_ic_id_i$ be the subpath of order $4$ of $T_2$ that contains $u_4$.

If either $u_4=b_i$ or $k=2$ and $u_4=c_i$, 
then $T\in {\cal T}$ and Lemma \ref{lemma4} implies the desired result.
Hence, we may assume that these cases do not occur.

Let $g_k$ be the function from the proof of Lemma \ref{lemma4}.
Clearly, $g_k$ is a Roman dominating function of $T$.
We will modifying $g_k$ on the four vertices in $\{ a_i,b_i,c_i,d_i\}$
in such a way that we obtain a Roman dominating function $g$ of $T$ of smaller weight than $g_k$.
Once this is done, we obtain $\gamma_{r2}(T)\leq \gamma_R(T)<3k=\frac{3}{4}n(T)$
and hence $\gamma_{r2}(T)+\gamma_R(T)<\frac{6}{4}n(T)$.

If $u_4=a_i$, then we change $g_k$ to $g$ such that $(g(a_i),g(b_i),g(c_i),g(d_i))=(0,0,2,0)$.
If $u_4=c_i$, then we may assume that $k\geq 3$, 
and we change $g_k$ to $g$ such that $(g(a_i),g(b_i),g(c_i),g(d_i))=(1,0,0,1)$.
Finally, 
if $u_4=d_i$, then we change $g_k$ to $g$ such that $(g(a_i),g(b_i),g(c_i),g(d_i))=(0,2,0,0)$.
It is straightforward to check that $g$ has the desired properties,
which completes the proof.
$\Box$

\medskip

\noindent Let the class ${\cal G}$ of connected graphs be such that 
a connected graph $G$ belongs to ${\cal G}$ if and only if $G$ arises
\begin{itemize}
\item either from the unique tree in ${\cal T}_2$ by adding the edge $c_1c_2$
\item or from some tree in ${\cal T}_k$ by arbitrarily adding edges between vertices in $\{ b_1,\ldots,b_k\}$,
\end{itemize}
where we denote the vertices as in the definition of ${\cal T}_k$.

\begin{theorem}\label{theorem2}
If $G$ is a connected graph of order $n(G)$ at least $3$, 
then $\gamma_{r2}(G)+\gamma_R(G)\leq \frac{6}{4}n(G)$
with equality if and only if $G\in {\cal G}$.
\end{theorem}
{\it Proof:} Similarly as in the proof of Lemma \ref{lemma4},
we obtain $\gamma_{r2}(G)+\gamma_R(G)=\frac{6}{4}n(G)$
for every graph $G\in {\cal G}$.
Now let $G$ be a connected graph of order at least $3$.
Let $T$ be a spanning tree of $G$.
By Lemma \ref{lemma5}, 
$\gamma_{r2}(G)+\gamma_R(G)\leq \gamma_{r2}(T)+\gamma_R(T)\leq \frac{6}{4}n(G)$.
Now, we assume that $\gamma_{r2}(G)+\gamma_R(G)=\frac{6}{4}n(G)$.
This implies $\gamma_{r2}(T)+\gamma_R(T)=\frac{6}{4}n(T)$.
By Lemma \ref{lemma5}, we obtain $T\in {\cal T}_k$ for some $k\in \mathbb{N}$.
We denote the vertices of $T$ as in the definition of ${\cal T}$.
Note that $T[\{ b_1,\ldots,b_k\}]$ is connected.
Let $f_k$ be the function from the proof of Lemma \ref{lemma4}.
Let $G_{\leq j}$ be the subgraph of $G$ induced by $\bigcup_{i=1}^j\{ a_i,b_i,c_i,d_i\}$.

First, we assume that $a_1$ has degree at least $2$ in $G$.
If $a_1$ is adjacent to $c_1$ or $d_1$, then we have $\gamma_{r2}(G_{\leq 1})=2$,
and replacing the values of $f_k$ on the vertices of $G_{\leq 1}$ 
by the values of an optimal $2$-rainbow dominating function of $G_{\leq 1}$
implies $\gamma_{r2}(G)<3k$, and hence 
$\gamma_{r2}(G)+\gamma_R(G)\leq \gamma_{r2}(T)+\gamma_R(T)<6k=\frac{6}{4}n(G)$,
which is a contradiction.
If $a_1$ is adjacent some vertex in $\{ a_2,b_2,c_2,d_2\}$,  
then we have $\gamma_{r2}(G_{\leq 2})\leq 5$,
and replacing the values of $f_k$ on the vertices of $G_{\leq 2}$ 
by the values of an optimal $2$-rainbow dominating function of $G_{\leq 2}$
implies $\gamma_{r2}(G)<3k$, which implies the same contradiction as above.
Hence, by symmetry, 
we may assume that all vertices in $\{ a_1,\ldots,a_k\}\cup \{ d_1,\ldots,d_k\}$ have degree $1$ in $G$.

If $b_1$ is adjacent to $c_2$, then we change the values of $f_k$ on the vertices in $\{ a_2,b_2,c_2,d_2\}$
to obtain a function $f$ with $(f(a_2),f(b_2),f(c_2),f(d_2))=(\{ 1\},\emptyset,\emptyset,\{ 1\})$.
Note that $b_2$ is adjacent to some vertex $b_i$ for $i\not=2$ and that $f_k(b_i)=\{ 1,2\}$.
Hence $f$ is a $2$-rainbow dominating function of $G$ and we obtain
$\gamma_{r2}(G)<3k$, which implies the same contradiction as above.
Hence, by symmetry, 
we may assume that no vertex $b_i$ is adjacent to some vertex $c_j$ for distinct $i$ and $j$.

If $b_1$ is adjacent to $b_2$ and $c_2$ is adjacent to $c_3$, 
then
$$
f:V(G_{\leq 3})\to 2^{\{ 1,2\}}:
u\mapsto
\left\{
\begin{array}{cl}
\emptyset & \mbox{, if $u\in \{ a_1,c_1,b_2,c_2,b_3,d_3\}$,}\\
\{ 1\} & \mbox{, if $u\in \{ d_1,a_2,d_2,a_3\}$, and}\\
\{ 1,2\} & \mbox{, if $u\in \{ b_1,c_3\}$}
\end{array}
\right.
$$
is a $2$-rainbow dominating function of $G_{\leq 3}$ of weight $8$.
This implies $\gamma_{r2}(G)<3k$, which implies the same contradiction as above.
Hence, by symmetry, 
we may assume that 
if some vertex $c_i$ has degree more than $2$,
then $b_i$ and $c_i$ both have degree exactly $3$ and 
there is some index $j$ distinct from $i$ such that 
$b_i$ is adjacent to $b_j$ 
and
$c_i$ is adjacent to $c_j$.

Altogether, the above observations easily imply that $G\in {\cal G}$ as desired.
$\Box$

\medskip

\noindent {\bf Acknowledgment}  
J.D. Alvarado and S. Dantas were partially supported by FAPERJ, CNPq, and CAPES.

\end{document}